\documentclass[11pt]{article}
\usepackage{amssymb}
\usepackage{amsmath}
\usepackage{amsthm}
\usepackage{latexsym}
\usepackage{amsfonts}
\usepackage{graphicx}
\usepackage{graphics}
\usepackage[T1]{pbsi}

\newtheorem{thm}{Theorem}
\newtheorem{lem}{Lemma}
\newtheorem{cor}{Corollary}
\newtheorem{Def}{Definition}

\theoremstyle{definition}
\newtheorem*{Proof}{Proof}

\newcommand{\dis}{\displaystyle}
\newcommand{\fa}{\forall}

\newcommand{\ra}{\;\rightarrow\;}

\newcommand{\al}{\alpha}

\newcommand{\Ga} {{\varGamma}}

\newcommand{\de}{\delta }

\newcommand{\f}{\varphi}

\newcommand{\la}{\lambda }
\newcommand{\mi}{\mu }

\newcommand{\R}{\mathbb{R}}

\newcommand{\N}{\mathbb{N}}

\newcommand{\ssum}{\sum\limits}

\newcommand{\ct}{{\cal{T}}}

\newcommand{\cm}{{\cal{M}}}

\newcommand{\ld}{\ldots}

\newcommand{\sm}{\smallsetminus}

\newcommand{\qb}{$\quad\blacksquare$}
 
 \newcommand{\ess}{\mbox{ess}}

\begin{document}

\title{\bf Dyadic $A_1$ weights and equimeasurable rearrangement of functions}
\author{E. Nikolidakis}
\date{}
\maketitle
\noindent
{\bf Abstract:} We prove that the decreasing rearrangement of a dyadic $A_1$-weight $w$ with dyadic $A_1$ constant $\big[w\big]^\ct_1=c$ with respect to a tree $\ct$ of homogeneity $k$, on a non-atomic probability space, is a usual $A_1$ weight on $(0,1]$ with $A_1$-constant $[w^\ast]_1$ not more\linebreak than $kc-k+1$. We prove also that the result is sharp, when one considers all such weights $w$. \vspace*{0.2cm} 
\section{Introduction}\label{sec1}
\noindent

The theory of Muckenhoupt weights has been proved to be an important tool in analysis due to their self-improving properties (see \cite{2,3,8}).

One class of special interest is $A_1(J,c)$ where $J$ is an interval on $\R$ and $c$ a constant $c\ge1$. Then $A_1(J,c)$ is defined as the class of all non-negative locally integrable functions $w$ defined on $J$, such that for every subinterval $I\subseteq J$ we have that
\begin{eqnarray}
\frac{1}{|I|}\int_Iw(y)dy\le c\,\ess\inf_{x\in I}w(x)  \label{eq1.1}
\end{eqnarray}
where $|\cdot|$ is the Lesbesgue measure on $\R$.

In \cite{1} it is proved that if $w\in A_1(J,c)$ then $w^\ast\in A_1((0,|J|],c)$, where $w^\ast$ is the non-increasing rearrangement of $w$. That is $w\in A_1(J,c)$ gives that
\begin{eqnarray}
\frac{1}{t}\int^t_0w^\ast(y)dy\le c\,w^\ast(t),  \label{eq1.2}
\end{eqnarray}
for every $t\in(0,|J|]$.

Here for a $w:J\ra\R^+$, $w^\ast$ stands for
\[
w^\ast(t)=\sup_{e\subseteq J\atop|e|\ge t}\inf_{x\in e}w(x), \ \ \text{for any} \ \ t\in (0,|J|].
\]
The fact mentioned above helps (as one can see in \cite{1}) in the determination of all $p$ such that $p>1$ and $w\in RH^J_p(c')$ for some $1\le c'<+\infty$ whenever $w\in A_1(J,c)$, where by $RH^J_p(c')$ we mean the class of all weights $w$ defined on $J$ which satisfy a reverse Holder inequality with constant $c'$ upon all the subintervals $I\subseteq J$. One can also see related problems for estimates for the range of $p$ in higher dimensions in \cite{4} and \cite{5}.

In this paper we are interested for the opposite dyadic case. A way of studying dyadic $A_1$ weights is by using the respective dyadic maximal operator.

More precisely, a locally integrable non-negative function $w$ on $\R^n$ is called a dyadic $A_1$ weight if it satisfies the following condition
\begin{eqnarray}
\frac{1}{|Q|}\int_Q w(y)dy\le c\,\ess\inf_{x\in I}w(x), \label{eq1.3}
\end{eqnarray}
for every dyadic cube on $\R^n$.

This condition is equivalent to the inequality
\begin{eqnarray}
\cm_d w(x)\le c\,w(x),  \label{eq1.4}
\end{eqnarray}
for almost all $x\in\R^n$. Here $\cm_d$ is the dyadic maximal operator defined by
\begin{eqnarray}
\cm_d w(x)=\sup\bigg\{\frac{1}{|Q|}\int_Q w(y)dy:x\in Q,\;Q\subset\R^n\;\text{is a dyadic cube}\bigg\}.  \label{eq1.5}
\end{eqnarray}
The smallest $c\ge1$ for which (\ref{eq1.3}) (equivalently (\ref{eq1.4})) holds is called the dyadic $A_1$ constant of $w$ and is denoted by $\big[w\big]^d_1$.

Let us now fix a dyadic cube $Q$ on $\R^n$. A natural problem that arises is the behaviour of $(w/Q)^\ast:(0,|Q|]\ra\R^+$ when one knows that $\big[w\big]^d_1=c$. It turns out that $(w/Q)^\ast$ is a usual $A_1$ weight on $(0,|Q|]$ with constant not more than $2^nc-2^n+1$.

More precisely we will prove the following
\begin{thm}\label{thm1}
Let $w$ be a dyadic $A_1$ weight on $\R^n$ with dyadic $A_1$ constant $\big[w\big]^d_1=c$. Let $Q$ be a fixed dyadic cube on $\R^n$. Then the following inequality is satisfied
\begin{eqnarray}
\frac{1}{t}\int^t_0(w/Q)^\ast(y)dy\le(2^nc-2^n+1)(w/Q)^\ast(t),  \label{eq1.6}
\end{eqnarray}
for every $t\in(0,|Q|]$.

Moreover the last inequality is sharp when one considers all dyadic $A_1$ weights with $\big[w\big]^d_1=c$.  \qb
\end{thm}

We remark that by using a standard dilation argument it suffices to prove (\ref{eq1.6}) for $Q=[0,1]^n$ and for all functions $w$ defined only on $[0,1]^n$ and satisfying the $A_1$ condition only for dyadic cubes contained in $[0,1]^n$. Actually, we will work on more general non-atomic probability spaces $(X,\mi)$ equipped with a structure $\ct$ similar to the dyadic one. (We give the precise definition in the next section).

The paper is organized as follows:

In Section \ref{sec2} we give some tools needed for the proof of Theorem \ref{thm1}. These are obtained from \cite{6} and \cite{7}.

In Section 3 we give the proof of Theorem \ref{thm1} in it's general form (as Theorem \ref{thm2}) and mention two applications of it.
\section{Preliminaries}\label{sec2}
\noindent

We fix a non-atomic probability space $(X,\mi)$ and a positive integer $k\ge2$.

We give the following
\begin{Def}\label{Def1}
A set of measurable subsets of $X$ will be called a tree of homogeneity $k$ if
\begin{enumerate}
\item[i)] For every $I\in\ct$ there corresponds a subset $C(I)\subseteq\ct$ containing exactly $k$ pairwise disjoint subsets of $I$ such that $I=\cup C(I)$ and each element of $C(I)$ has measure $(1/k)\mi(I)$.
\item[ii)] $\ct=\bigcup\limits_{m\ge0}\ct_{(m)}$ where $\ct_{(0)}=\{X\}$ and $\ct_{(m+1)}=\bigcup\limits_{I\in\ct_{(m)}}C(I)$.
\item[iii)] The tree $\ct$ differentiates $L^1(X,\mi)$, that is if $\f\in L^1(X,\mi)$ then for $\mi$ almost all $x\in X$ and every sequence $(I_k)_{k\in\N}$ such that $x\in I_k$, $I_k\in\ct$ and $\mi(I_k)\ra0$ we have that
\[
\f(x)=\lim_{k\ra+\infty}\frac{1}{\mi(I_k)}\int_{I_k}\f d\mi. \text{\qb}
\]
\end{enumerate}
\end{Def}

It is clear that each family $\ct_{(m)}$ consists of $k^m$ pairwise disjoint sets, each having measure $k^{-m}$, whose union is $X$.

Moreover, if $I,J\in\ct$ and $I\cap J$ is non empty then $I\subseteq J$ or $J\subseteq I$.

For this family $\ct$ we define the associated maximal operator $\cm_\ct$ by
\setcounter{equation}{0}
\begin{eqnarray}
\cm_\ct\f(x)=\sup\bigg\{\frac{1}{\mi(I)}\int_I|\f|d\mi:\,x\in I\in\ct\bigg\}, \label{eq2.1}
\end{eqnarray}
for any $\f\in L^1(X,\mi)$ and we will say that a non-negative integrable function $w$ is an $A_1$ weight with respect to $\ct$ if
\begin{eqnarray}
\cm_\ct\f(x)\le C\f(x),  \label{eq2.2}
\end{eqnarray}
for almost every $x\in X$. The smallest constant $C$ for which (\ref{eq2.2}) holds will be called the $A_1$ constant of $w$ with respect to $\ct$ and will be denoted by $\big[w\big]^\ct_1$.

We give now the following:
\begin{Def}\label{Def2}
Every non-constant function $w$ of the form $w=\ssum_{P\in\ct_{(m)}}\la_P\xi_P$, for a specific $m>0$, and for positive $\la_P$, will be called a $\ct$-step function. ($\xi_P$ denotes the characteristic function of $P$). \qb
\end{Def}

It is then clear that every $\ct$-step function is an $A_1$ weight with respect to $\ct$. Let $\de=1/\big[w\big]^\ct_1$, $0<\de<1$ and for any $I\in\ct$ write $Av_I(w)=\frac{1}{\mi(I)}\int\limits_Iwd\mi$.

Now for every $x\in X$, let $I_w(x)$ denote the largest element of the set $\{I\in\ct:\,x\in I$ and $\cm_\ct w(x)=Av_I(w)\}$ (which is non-empty since $Av_J(w)=Av_P(w)$ for every $P\in\ct_{(m)}$ and $J\subseteq P$).

Next for any $I\in\ct$ we define the set
\[
A_I=A(w,I)=\{x\in X:\,I_w(x)=I\}
\]
and let $S=S_w$ denote the set of all $I\in\ct$ such that $A_I$ is non-empty. It is clear that each such $A_I$ is a union of certain $P$ from $\ct_{(m)}$ and moreover
\[
\cm_\ct w=\sum_{I\in S}Av_I(w)\xi_{A_I}.
\]
We also define the correspondence $I\ra I^\ast$ with respect to $S$ as follows: $I^\ast$ is the smallest element of $\{J\in S_w:\,I\varsubsetneq J\}$. This is defined for every $I\in S$ that is not maximal with respect to $\subseteq$.

We recall parts of two Lemmas from \cite{6}.
\begin{lem}\label{lem1}
For all $I\in\ct$ we have $I\in S$, if and only if, $Av_Q(w)<Av_I(w)$ whenever $I\subseteq Q\in\ct$, $I\neq Q$. In particular $X\in S$ and so $I\ra I^\ast$ is defined for all $I\in S$ such that $I\neq X$.  \qb
\end{lem}
\begin{lem}\label{lem2}
Let $I\in S$. Then, if $J\in S$ is such that
\[
J^\ast=I \ \ \text{then} \ \ y_I<y_J\le(k-(k-1)\de)y_I.  \text{\qb}
\]
\end{lem}
\section{Main theorem and proof}\label{sec3}
\noindent

In this section we will prove the following.
\begin{thm}\label{thm2}
Let $\ct$ be a tree of homogeneity $k\ge2$ on the probability non-atomic space $(X,\mi)$, and let $w$ be $A_1$ weight with respect to $\ct$ with $A_1$-constant $\big[w\big]^\ct_1=c$. Then if one considers $w^\ast:(0,1]\ra\R^+$ we have that $\frac{1}{t}\int\limits_0^tw^\ast(y)dy\le(kc-k+1)w^\ast(t)$, for every $t\in(0,1]$, where as usual $w^\ast$ is defined by $w^\ast(t)=\dis\sup_{e\subseteq x\atop\mi(e)\ge t}\dis\inf_{x\in e}w(x)$, $t\in(0,1]$.

Moreover the constant appearing in the right of the last inequality is sharp, if one considers all $A_1$ weights with respect to $\ct$ with constant $\big[w\big]^\ct_1=c$. \qb
\end{thm}
\begin{Proof}
We suppose for the beginning that $w$ is a $\ct$-step function. Fix $t\in(0,1]$ and consider the set
\begin{align*}
E_t&=\{x\in X:\,\cm_\ct w(x)>c\,w^\ast(t)\} \\
&=\{\cm_\ct w>c\la\}, \ \ \text{where} \ \ \la=w^\ast(t).
\end{align*}
Then $E_t$ is a measurable subset of $X$. We first assume that $\mi(E_t)>0$.

We consider the family of all those $I\in\ct$ maximal under the condition $Av_I(w)>c\la$, and denote it by $(I_j)_j$. Then $(I_j)_j$ is pairwise disjoint and $E_t=\cup I_j$.

Additionally for every $j$ and $I\in\ct$ such that $I\supsetneq I_j$ we have that $\frac{1}{\mi(I)}\int\limits_Iwd\mi=Av_I(w)\le c\la$ because of the maximality of $I_j$.

In view of Lemma \ref{lem1} this gives $I_j\in S_w=S$, for every $j$.

For every $I_j$ consider $I^\ast_j\in S$. Then by Lemma \ref{lem2}, $y_{I_j}\le[k-(k-1)\de]y_{I^\ast_j}$, where $\de=1/c$ and of course $y_{I^\ast_j}\le c\la$. So, we have that
\[
y_{I_j}\le[k-(k-1)\de]c\la=(kc-k+1)\la, \ \ \text{for every} \ \ j.
\]
This gives
\setcounter{equation}{0}
\begin{align}
\int_{I_j}wd\mi\le(kc-k+1)\la\mi(I_j)&\Rightarrow\int_{E_t}wd\mi\le(kc-k+1)\la\mi(E_t)\nonumber\\
&\Rightarrow\frac{1}{\mi(E_t)}\int_{E_t}wd\mi\le(kc-k+1)\la. \label{eq3.1}
\end{align}
Since $\cm_\ct w\le cw$ on $X$, and $E_t=\{\cm_\ct w>c\la\}$ we obviously have $E_t\subseteq\{w>\la\}=\{w>w^\ast(t)\}$.

There exist now $E^\ast_t\subseteq(0,1]$ Lesbesgue measurable such that $|E^\ast_t|=\mi(E_t)=:t_1$, and such that $\int\limits_{E^\ast_t}w^\ast(y)dy=\int\limits_{E_t}wd\mi$. Obviously we can arrange everything in a way such that $E^\ast_t\subseteq\{w^\ast>w^\ast(t)\}\subseteq(0,t)$. As a result $t_1\le t$.

Since now $\ct$ differentiates $L^1(X,\mi)$ we have that almost every element of the set $\{w>c\la\}\subseteq X$ belongs to $E_t$.

Since $\mi(E_t)>0$ we have that $\mi(\{w>c\la\})>0$.

Let now $t_2$ be such that
\[
w^\ast(t)>\la c \ \ \text{for every} \ \ t\in (0,t_2) \ \ \text{and} \ \ w^\ast(t)\le c\la, \ \ \text{for every} \ \ t\in(t_2,1).
\]
Then, we can arrange everything (by deleting suitable sets of Lesbesgue measure zero) in a way that $E^\ast_t=(0,t_2)\cup A_t$, where $A_t$ is a Lesbesgue measurable subset of $(t_2,t)$ and $|A_t|=t_1-t_2$ (Of course $t_2=|(0,t_2)|=|\{w^\ast>\la c\}|=\mi(\{w>\la c\})\le\mi(\{\cm_\ct w>\la c\})=\mi(E_t)=:t_1)$.

We will prove the following
\begin{eqnarray}
\frac{1}{\mi(E_t)}\int_{E_t}wd\mi\ge\frac{1}{t}\int^t_0w^\ast(y)dy,  \label{eq3.2}
\end{eqnarray}
(\ref{eq3.2}) is equivalent to
\begin{align}
\frac{1}{t_1}\int_{E^\ast_t}w^\ast(y)dy\ge&\,\frac{1}{t}\int^t_0w^\ast(y)dy\Leftrightarrow t
\int^{t_2}_0w^\ast(y)dy+t\int_{A_t}w^\ast(y)dy\nonumber \\
\ge&\, t_1\int^{t_2}_0w^\ast(y)dy+t_1\int^t_{t_2}w^\ast(y)dy \nonumber\\
&\Leftrightarrow(t-t_1)\int^{t_2}_0w^\ast(y)dy+t\int_{A_t}w^\ast(y)dy\nonumber \\
\ge&\, t_1\int^t_{t_2}w^\ast(y)dy,  \label{eq3.3}
\end{align}
We define $\Ga_t=(t_2,t)\sm A_t$. (\ref{eq3.3}) then becomes
\begin{align}
&(t-t_1)\int^{t_2}_0w^\ast(y)dy+(t-t_1)\int_{A_t}w^\ast(y)dy\ge t_1\int_{\Ga_t}w^\ast(y)dy \nonumber \\
&\Leftrightarrow(t-t_1)\int_{E^\ast_t}w^\ast(y)dy\ge t_1\int_{\Ga_t}w^\ast(y)dy.  \label{eq3.4}
\end{align}
But of course
\[
\int_{E^\ast_t}w^\ast(y)dy=\int_{E_t}wd\mi>\mi(E_t)\cdot c\la=c\la\cdot t_1,
\]
in view of the known weak type inequality for $\cm_\ct$, namely:
\[
\mi(\{\cm_\ct\f>\la\})<\frac{1}{\la}\int_{\{\cm_\ct\f>\la\}}\f.
\]
So, if we prove that
\begin{eqnarray}
\int_{\Ga_t}w^\ast(y)dy\le c\la(t-t_1),  \label{eq3.5}
\end{eqnarray}
we complete the proof of (\ref{eq3.2}). But (\ref{eq3.5}) is obvious since $w^\ast(y)\le c\la$ on $(t_2,t)$, $\Ga_t\subseteq(t_2,t)$ and
\[
|\Ga_t|=|(t_2,t)\sm(A_t)|=(t-t_2)-|A_t|=t-t_1.
\]
We have thus proved for every $w$ $\ct$-step function and $t$ such that $\mi(\{\cm_\ct w>c\cdot w^\ast(t)\})>0$, that
\begin{eqnarray}
\frac{1}{t}\int_0^tw^\ast(y)dy\le(kc-k+1)w^\ast(t).  \label{eq3.6}
\end{eqnarray}
If $t$ is such that $\mi(\{\cm_\ct w>cw^\ast(t)\})=0$ then obviously $\cm_\ct w_{(x)}\le cw^\ast(t)$, for almost all $x\in X$, so since $\ct$ differentiates $L^1(X,\mi)$: $w(y)\le cw^\ast(t)$ for almost all $y\in X$. This obviously give (\ref{eq3.6}) since $c\le kc-k+1$.

Additionally if $w$ is in general an $A_1$-weight with respect to $\ct$, then an approximation argument by $\ct$-simple $A_1$-weights gives the result for $w$.

More precisely one can easily see, that if $w$ is a $A_1$ weight with respect to $\ct$, with $A_1$-constant $\big[w\big]^\ct_1=c$ then there exist a sequence of $\ct$-simple functions $(w_n)_n$ increasing as $n$ increases, and such that $w_n\le w$ and $\big[w\big]^\ct_1=c_n\le c$ with $w_n\ra w$ and $c_n\ra c$ as $n\ra+\infty$.

In order to finish the proof of Theorem \ref{thm2} we just need to prove the sharpness of the result. We do it right now:

Fix $k\ge2$. We suppose that we are given a tree $\ct$ of homogeneity $k$, and consider $\ct_{(2)}$. Then
\[
\ct_{(2)}=\{P_1,\ld,P_k,P_{k+1},\ld,P_{2k},\ld,P_{k^2-k+1},\ld,P_{k^2}\} \ \ \text{where}
\]
\[
\ct_{(1)}=\bigg\{\bigcup^k_{i=1}P_i,\bigcup^{2k}_{i=k+1}P_i,\ld,\bigcup^{k^2}_
{i=k^2-k+1}P_i\bigg\}=\{I_{1},I_2,\ld,I_k\}.
\]
We have that $\mi(P_i)=\frac{1}{k^2}$, $\fa\;i$.

Suppose $\de>0$ be such that $\de<\frac{1}{k^2}$, and consider for any such $\de$ a set $A_\de$ of measure $\mi(A_\de)=\de$ such that $A_\de\subseteq P_1$ ($(X,\mi)$ is non atomic). Let $c\ge1$ and $\al,\epsilon<0$. Let $\f=\f_\de$ be the function defined as follows:
\[
\begin{array}{ll}
  \f/A_\de:=\al &  \\
  \f/I_1\sm A_\de:=\epsilon &  \\
  \f/P_{k+1}:=\al, & \f/I_2\sm P_{k+1}:=\epsilon \\
   \f/P_{2k+1}:=\al, & \f/I_3\sm P_{2k+1}:=\epsilon\\
\multicolumn{2}{c}{\cdots}   \\
  \f/P_{k^2-k+1}:=\al, & \f/I_k\sm P_{k^2-k+1}:=\epsilon \\
\end{array}
\]
It is easy to see that $\f=\f_\de$ is a $A_1$ weight with $A_1$ constant %
\[
c_\de=\big[w\big]^\ct_1=\frac{Av_{I_1}(\f)}{\epsilon}=\frac{k}{\epsilon}\int_I
\f d\mi=\frac{k}{\epsilon}\bigg[a\de+\bigg(\frac{1}{k}-\de\bigg)\epsilon\bigg].
\]
Then $c_\de\ra c$, as $\de\ra1/k^{2-}$ iff: $\al,\epsilon$ are chosen such that $kc-k+1=\frac{\al}{\epsilon}$. (Given $k,c$). Let us choose $\al,\epsilon$ be such as mentioned just before, with $\epsilon<\al$.

Then $\f^\ast_\de(1/k)=\epsilon$, so $\f^\ast_\de(1/k)(kc-k+1)=\al$, while $k\int^{1/k}_0\f^\ast_\de(y)dy$ tends to $\al$, while $\de\ra1/k^{2^-}$.

By this we end the proof of Theorem \ref{thm2}.  \qb
\end{Proof}

Theorem \ref{thm1} of Section 1 is an immediate Corollary of Theorem \ref{thm2}.

Additionally the following are consequences of Theorem \ref{thm2}.
\begin{cor}\label{cor1}
Let $w$ be an $A_1$ weight with respect to the tree $\ct$ of homogeneity $(k\ge2)$ on $(X,\mi)$ with $\big[w\big]^\ct_1=c$. Then if one considers $((0,1],|\cdot|)$ equipped with the usual $k$-adic tree $\ct_k$, where $|\cdot|$ is the Lesbesgue measure on $(0,1]$. Then $\big[w^\ast\big]^{\ct_k}_1\le kc-k+1$ and this result is sharp.
\end{cor}
\begin{Proof}
Obvious, according to the function $\f_\de$ constructed at the end of Theo-\linebreak rem \ref{thm2}. \qb
\end{Proof}
\begin{cor}\label{cor2}
Let $w$ be $A_1$-weight on $\R^n$ as described in Section \ref{sec1} Then $w^\ast:(0,+\infty)\ra\R^+$ has the following property:
\[
\frac{1}{t}\int^t_0w^\ast(y)dy\le(kc-k+1)w^\ast(t), \ \ \text{for every} \ \ t\in(0,+\infty)
\]
and the last inequality is sharp.  \qb
\end{cor}
\begin{Proof}
We expand $\R^n$ as a union of an increasing sequence $(Q_j)_j$ of dyadic cubes, and use Theorem \ref{thm2} in any of these. \qb
\end{Proof}
\bigskip\medskip
Eleftherios Nikolidakis  \vspace*{0.1cm} \\
Department of Mathematics \vspace*{0.1cm} \\
University of Crete,\vspace*{0.1cm} \\
Heraclion 71409,  \vspace*{0.1cm} \\
Crete, Greece\vspace*{0.1cm} \\
E-mail address: lefteris@math.uoc.gr

\end{document}